\documentclass[12pt]{article}
\usepackage{amsmath}
\usepackage{amssymb}

\usepackage[cp1251]{inputenc}
\usepackage[russian]{babel}
\newcommand{\il}[2]{\int\limits_{#1}^{#2}}
\newcommand{\ilp}[1]{\int\limits_{#1}^{+\infty}}

\newcommand{\ph}{\phantom{a}}
\newcommand{\phh}{\phantom{aaa}}

\newcommand{\sist}[2]{\left\{
\begin{array}{l}
{#1}\\
\ph\\
{#2}
\end{array}
\right.}

\begin{document}

MSC 47B35

\vskip 15pt

\centerline{\bf Solvabilty conditions for a special class}
\centerline{\bf   of Wiener-Hopf  integral   equations}

\vskip 10pt
\centerline{\bf G. A. Grigorian}

\centerline{\it Institute  of Mathematics of NAS of Armenia}
\centerline{\it E -mail: mathphys2@instmath.sci.am}
\vskip 20 pt

Abstract. Inthis paper we study the Wiener-Hopf integral equation in the case of kerkel function hanving a special form of mixture of Gaussian, Cauchy and Poisson superpositions, havimg importanse in the probability theory. We apply the classical theory for elliptic Wiener-Hopf equations, as well as the method of special factoization to obtain solvability conditions for Eq. (1.1)

\vskip 10pt


\vskip 10pt

{\bf 1. Introduction}. Let $K(t)$ be an absolutely integrable function on $\mathbb{R}$, and let $f(t)$ be a function on $[0,\infty)$. Consider the Wiener-Hopf  integral equation
$$
\mu \phi(t) -  \int\limits_{0}^\infty K(t - \tau) \phi(\tau) d \tau = f(t), \ph t \ge 0,  \ph \mu\in\mathbb{R}.
 \eqno (1.1)
$$
Without loss of generality we will assume that
$$
\int\limits_{-\infty}^\infty K(t)dt=1. \eqno (1.2)
$$
This equation appears in the study of many problems of engineering and physics,
  In the case      $K(t)=K(-t), \ph  -\infty< t<\infty,  \ph 0<\mu\le1$  it relates to the class of Wiener-Hopf equations of non normal type, since its symbol degenerates. In [1] Eq. (1.1) was studied in the case $0<\mu<1$ under the additional restriction $0\le K\in L_1(\mathbb{R}^+\cap C^{(2)}(\mathbb{R}^+), \linebreak  K'\le ~0, \ph 0\le K''\downarrow$ on $\mathbb{R}^+$. In [2] Eq. (1.1) was studied in the case $0<\mu<1$ with the additional assumption that $\mu$ is enough close to $1$.
These restrictions are too stringent to reveal the richness of the qualitative diversity of the properties of equation (1.1) regarding the existence, uniqueness, and asymptotic behavior of its solutions.

In this paper we study the solvability conditions of Eq. (1.1)  in the case, when the kernel function $K(t)$ has the following special form
$$
K(t)=\sist{\int\limits_{0}^\infty e^{-st^2} d\sigma_{1}^+(s) + \int\limits_{0}^\infty\frac{d\sigma_{2}^+(s)}{t^2 + s^2} + \int\limits_{0}^\infty e^{-st} d\sigma_{3}^+(s), \ph t \ge 0,}{\int\limits_{0}^\infty e^{-st^2} d\sigma_{j}^-(s) +
\int\limits_{0}^\infty\frac{d\sigma_{2}^-(s)}{t^2 + s^2} +  \int\limits_{0}^\infty e^{st} d\sigma_{3}^-(s),, \ph t< 0, }  \eqno (1.3)
$$
where
$\sigma_j(t), \ph j=\overline{1,3}$ are real-valued  functions of local finite variation  on $[0,\infty)$, satisfying the conditions

\noindent
1) $\int\limits_{0}^\infty\frac{|d\sigma_1^\pm(s)|}{\sqrt{s}} < \infty,\phh$
2) $\int\limits_{0}^\infty\frac{|d\sigma_j^\pm(s)|}{s} d s < \infty, \ph j=2,3.$

\vskip 10pt

Note that  $K(t)$ is a mixture of compositions of  Gaussian, Cauchy and Poisson's  distributions, playing important role in the probability theory.
The special form (1.3) of $K(t)$, as it is shown below, allowed us to overcome restrictions on $K(t)$,imposed on it in [1] and [2]. We apply the classical theory of elliptic Wiener-Hopf equations, as well as the method of special factorization to obtain solvability conditions for Eq. (1.1).  This equation has important applications in probability theory and astrophysics   in the  so called ``conservative'' case ($\mu=1$), to which devoted many investigations (see e.g. [3, 4] and cited works  therein).
This case was studied for the general kernel function   $K(t)$  with $0\le K\in L_1(\mathbb{R})$. The special structure (1.3) of $K(t)$  allowed us to obtain qualitatively new results, (different from results of works [3,4] and others, see below Theorems  3.4, 3.5, where the last restriction on $K$ may not be  fulfilled, namely  $K(t)+K(-t), \ph t\in \mathbb{R}$ may change sign   [see below Examples 2.1, 3.1]).

\vskip 10pt

{\bf 2. Auxiliary propositions}.

 For any function    $k\in L^1(\mathbb{R})$     denote by $\widehat{K}(\lambda)$ its Fourier transform
$$
\widehat{K}(\lambda) \equiv \int\limits_{-\infty}^{\infty} e^{i\lambda t} k(t) d t, \phh \lambda \in \mathbb{R}.
$$
By $\mathbb{W}$ we denote the extended Wiener algebra of Fourier transforms of absolutely integrable functions $k(t)$ on $\mathbb{R}$: $\mathbb{W} \equiv \{c+ \widehat{K}(\lambda), \ph c\in \mathbb{C}\}$.

 For any $k\in L^1(\mathbb{R})$ we set \ph
$
\nu_m(k) \equiv \il{-\infty}{\infty} t^m k(t) d t, \ph m= 0, 1, 2, \dots.
$
Denote by $E$ one of the spaces $L^p(\mathbb{R}) \ph (1 \le p < \infty), \ph \mathbb{C}^0 \subset \mathbb{M}$, where $\mathbb{M}$ is the space of measurable essentially bounded on $\mathbb{R}$ functions, $\mathbb{C}^0$ is the space of continuous on $\mathbb{R}$ functions $g(t)$ with $\lim\limits_{t \to \pm\infty} g(t) = 0$. We denote by $E_+$ the restriction of $E$ on $[0,\infty)$, i. e.
$$
E_+ = \{P\phi, \ph \phi \in E\},
$$
where $P$ is the following orthogonal projector:
$$
(P\phi)(t) = \sist{\phi(t), \ph t > 0,}{0, \ph t <0.}
$$
Let $\alpha$ be a real number. We set:
$$
(B_\alpha \phi)(t) \equiv \phi(t) - (1+ i \alpha) e^t\il{t}{\infty}e^{-s}\phi(s) d s, \ph \phi \in E_+, \ph t \ge 0,
$$
$$
(B_\infty \phi)(t) \equiv i e^t \il{t}{\infty} e^{- s} \phi(s) d s, \ph \phi \in E_+, \ph t \ge 0,
$$
$$
(G_\alpha \phi)(t)\equiv \phi(t) + (1 - i\alpha) e^{-i \alpha t} \il{0}{t} e^{i \alpha s} \phi(s) d s, \phi \in E_+, \ph t \ge 0,
$$
$$
(G_\infty f)(t) \equiv i [f(t) + f'(t)], \ph t \ge t_0,
$$
where $f'(t)$ is the generalized derivative of $f(t)$.
Let $\alpha_1, \ldots, \alpha_r$ be any real numbers,  $m_1, \ldots, m_r, m_\infty$ be any natural numbers   (not all necessarily different),    and let $m_1', \ldots, m_r', \linebreak m_\infty, \ph m_1'', \ldots, m_r'', \ph m_\infty''$  be nonnegative integers, such that
$$
m_1' + m_1'' = m_1, \ldots m_r' + m_r'' = m_r, \ph m_\infty' + m_\infty'' = m_\infty.
$$

We denote
$$
B\equiv B_\infty^{m_\infty''}\prod_{j=1}^r B_{\alpha_j}^{m_\infty''}, \phh G\equiv G_\infty^{m_\infty''}\prod_{j=1}^r G_{\alpha_j}^{m_\infty''},
$$
$$
\rho_+(\lambda) \equiv (\lambda + i)^{-m_\infty'}\prod\limits_{j=1}^r \biggl(\frac{\lambda - \alpha_j}{\lambda + i}\biggr)^{m_j'}, \ph \rho_-(\lambda) \equiv (\lambda - i)^{-m_\infty''}\prod\limits_{j=1}^r \biggl(\frac{\lambda - \alpha_j}{\lambda - i}\biggr)^{m_j''},
$$
$-\infty \le \lambda \le \infty$. The functions $\rho_+(\lambda)$ and $\rho_-(\lambda)$ are called the symbols of the operators $G$ and $B$ respectively (see [6, p. 175]) and the spaces
$$
\widetilde{E}_+(\rho_+) \equiv G(E_+), \phh \overline{E}_+(\rho_-) \equiv B(E_+)
$$
are called the spaces, generated by the zeroes of the symbols  $\rho_+$ and $\rho_-$ respectively. The norm $||f||_{\overline{E}_-}$ in $\overline{E}_+(\rho_-)$ is defined as follows
$$
||f||_{\overline{E}_-} = ||B^{-1} f||_{E_+} \ph \mbox{for} \ph  E_+ = L^P_+ (1 \le p < \infty) \ph \mbox{or} \ph E_+ = C_+^0
$$
and
$$
||f||_{\overline{E}_-} = \inf\limits_{B\phi = f} || \phi||_{E_+}, \ph \mbox{for} \ph E_+ = M_+.
$$
By similar way it is defined the norm $||f||_{\widetilde{E}_+}$ in $\widetilde{E}_+(\rho_+)$.

\vskip 10pt

 Consider the Wiener-Hopf operator
$$
(\mathbb{A}\phi) \equiv c \phi(t) - \int\limits_{0}^\infty k(t-s) \phi(s) d s, \phh t \ge 0
$$
and the corresponding symbol
$$
\widehat{A}(\lambda) \equiv c - \il{-\infty}{\infty} e^{ i \lambda t} k(t) d t \ph (-\infty \le \lambda \le \infty).
$$
Assume
$$
\widehat{A}(\lambda) = \rho_+(\lambda) \mathcal{C}(\lambda) \rho_-(\lambda), \phh \mathcal{C}(\lambda) \in \mathbb{W},  \phh -\infty \le \lambda \le  \infty.
$$

\vskip 10pt

{\bf Theorem 2.1} ([5, p. 175, Theorem 2.6]). {\it  For the operator $\mathbb{A}$ to be a $\Phi_+$-operator or a $\Phi_-$-operator from $\widetilde{E}_+(\rho_+)$ into $\overline{E}(\rho_-)$ it is necessary and sufficient that $\mathcal{C}(\lambda) \ne 0, \ph -\infty \le \lambda \le \infty.$ If this condition is satisfied, then the invertibility of the operator $\widehat{A}$ corresponds to the number
$\kappa - \delta$, and the formulae
$$
dimker \mathbb{A} = \max \{\delta -\kappa, 0\}, \phh dimcoker \mathbb{A} = \max\{\kappa - \delta, 0\}
$$
hold. Here $\kappa = ind \mathcal{C}(\lambda),$
$$
\delta = \sist{0, \ph \mbox{for} \ph E_+ = L_+^p \ph (1 \le p <\infty) \ph \mbox{or} \ph E_+ = C_+^0,}{r'', \ph \mbox{for} \ph E_+ = M_+,}
$$
and $r''$ is the number of positive $m_j'' (j=1,2,\ldots,r).$ In the case $E_+ = C_+$ we have
$$
\delta = \sist{1, \ph \mbox{if} \ph \alpha_{j_0} = 0 \ph \mbox{and} \ph m_{j_0}'' > (1 \le j_0 \le r),}{0, \ph \mbox{otherwise}.}
$$
}

\phantom{aaaaaaaaaaaaaaaaaaaaaaaaaaaaaaaaaaaaaaaaaaaaaaaaaaaaaaaaaaaaaaaaaaaaaa} $\blacksquare$

Let us introduce an operator $R_\alpha (-\infty < \alpha < \infty)$, which is defined by the formula
$$
(R_\alpha\phi)(t)=\sist{i e^{-i\alpha t}\int\limits_{t}^\infty e^{i\alpha s } \phi(s) ds, \ph t > 0,}{i e^{-i\alpha t}\int\limits_{
-\infty}^t e^{i\alpha s} \phi(s) ds, \ph t <0,}
$$
$\phi\in L^1(-\infty,\infty)$. Obviously,
$$
if \ph \nu_s(|\phi|)<\infty, \ph, s=0,1,\ldots m, \ph then \ph (R_\alpha^s\phi)(t)\in L^1,  \ph s=1,2,\ldots, m. \eqno (2.1)
$$

\vskip 10pt

Let $\alpha_0(=0), \ph \alpha_1, \ldots, \alpha_n$ be different  roots of $\widehat{A}(\lambda)$ and  $m_0, \ph m_1, \ldots m_n$
be any nonnega-\linebreak tive integers.

\vskip 10pt

{\bf Theorem 2.2([5, p. 180, Theorem 2.9]).} {\it In order that
$$
\widehat{A}(\lambda)=\rho(\lambda)\mathcal{B}(\lambda),
\ph \mathcal{B}(\lambda) \in\mathbb{W},
$$
where
$$
\rho(\lambda) \equiv \Bigl(\frac{\lambda}{\lambda + i}\Bigr)^{m_0}\prod\limits_{j=1}^n\Bigl(\frac{\lambda - \alpha_j}{\lambda + i }\Bigr)^{m_j},
$$
it is necessary and sufficient that
$$
(R_{\alpha_j}^s k)(t)\in L^1, \ph s=0,1,\ldots m_j, \ph j=0,1,\ldots n.
$$
}

\phantom{aaaaaaaaaaaaaaaaaaaaaaaaaaaaaaaaaaaaaaaaaaaaaaaaaaaaaaaaaaaaaaaaaaaaaaaa} $\blacksquare$




We set
$$
{K}_1(t) \equiv\sist{\int\limits_{0}^\infty e^{-st^2} d \sigma_1^+(s), \ph t \ge 0,}{\int\limits_{0}^\infty e^{-st^2} d \sigma_1^-(s), \ph t< 0,} \phh {K}_2(t) \equiv\sist{\int\limits_{0}^\infty \frac{d \sigma_2^+(s)}{t^2 + s^2}, \ph t \ge 0,}{\ilp{0}\frac{d \sigma_2^-(s)}{t^2 + s^2}, \ph t < 0,}
$$
$$
{K}_3(t) \equiv\sist{\int\limits_{0}^\infty e^{-st} d\sigma_3^+(s), \ph t \ge 0,}{\int\limits_{0}^\infty e^{st} d\sigma_3^-(s), \ph t< 0,}
$$
Then, obviously,
$$
K(t) = K_1(t) + K_2(t) + K_3(t), \ph t \in \mathbb{R}. \eqno (2.2)
$$

On the basis of the well known equality for the Poisson's integral (see e.g., [6, p. 415])
$$
\int\limits_{0}^\infty e^{-ax^2}\cos bx d x = \frac{1}{2}\sqrt{\frac{\pi}{a}} e^{-b^2/4a} \ph (a > 0)
$$
it is easy to derive
$$
\int\limits_{-\infty}^{\infty}e^{i\lambda t} e ^{-s t^2} d t = \sqrt{\frac{\pi}{s}} e^{-\lambda^2/4s} \ph (s > 0)
$$
Using this equality one can show that
$$
\widehat{{K}}_1(\lambda) = \int\limits_{0}^\infty\sqrt{\frac{\pi}{s}}  e^{-\lambda^2/4s}d\Bigl( \frac{\sigma_1^+(s) + \sigma_1^-(s)}{2}\Bigr) + i \int\limits_{0}^{\infty}\sin \lambda t   d t    \int\limits_{0}^\infty e^{-s t^2} d[\sigma_1^+(s) - \sigma_1^-(s)],  \eqno  (2.3)
$$
$ \lambda \in \mathbb{R}$.
Since $\int\limits_{0}^\infty e^{i\lambda t  - st} d t = \frac{-1}{i \lambda - s} \ph,  \il{-\infty}{0}e^{i\lambda t  + st} d t = \frac{1}{i \lambda + s}$ we have
$$
\int\limits_{-\infty}^{\infty}e^{i \lambda t}\frac{e^{-s|t|}}{2 s} d t = \frac{1}{s^2 + \lambda^2} \phh (s > 0).
$$
By the formula for inverse Furier transformation from here we obtain
$$
\int\limits_{-\infty}^{\infty}\frac{e^{i\lambda t}}{s^2 + t^2} d t = \frac{e^{-s|\lambda|}}{4\pi s}, \phh \lambda \in \mathbb{R}, \ph s > 0.
$$
Using this equality one can show that
$$
\widehat{{K}}_2(\lambda) =
\int\limits_{0}^\infty  e^{-s|\lambda|}\frac{d(\sigma_2^+(s) + \sigma_2^-(s))}{8\pi s} + i \int\limits_{0}^\infty\sin \lambda t d t\int\limits_{0}^\infty\frac{d(\sigma_2^+(s) - \sigma_2^-(s))}{2(t^2 + s^2)}, \phh \lambda \in \mathbb{R}. \eqno (2.4)
$$
Using the equalities $\int\limits_{0}^\infty e^{i\lambda t  - st} d t = \frac{1}{s-i \lambda} \ph,  \il{-\infty}{0}e^{i\lambda t  + st} d t = \frac{1}{i \lambda + s} (s > 0)$ one can also show that
$$
\widehat{{K}}_3(\lambda) = \int\limits_{0}^\infty\frac{s}{s^2  + \lambda^2} d(\sigma_3^+(s) + \sigma_3^-(s))   + i \lambda \int\limits_{0}^\infty\frac{s}{s^2 + \lambda^2}d(\sigma_3^+(s) - \sigma_3^-(s)), \ph \lambda \in \mathbb{R}, \eqno (2.5)
$$



Before starting the study the symbol of Eq. (1.1) let us prove the following lemma.

\vskip 10pt

{\bf Lemma 2.1.} {\it Let the following conditions be satisfied.
$$
A)
\int\limits_{t}^\infty\frac{d(\sigma_1^+(s) + \sigma_1^-(s))}{s^{3/2}} \ge 0, \phh
B) \il{0}{t}d(\sigma_j^+(s) + \sigma_j^-(s)) \ge 0, \ph t \ge 0, \ph j=2,3,
$$

 Then $Re \hskip 2pt \widehat{{K}}(\lambda)$ nondecreases on $(-\infty,0]$, nonincreases  on $[0,\infty)$ and
$$
\max\limits_{\lambda \in \mathbb{R}} Re \hskip 2pt \widehat{{K}}(\lambda) =\widehat{{K}}(0) = \nu_0({K}).
$$
Moreover,  if
 $$
C) \int\limits_{t}^\infty\frac{d(\sigma_1^+(s) + \sigma_1^-(s))}{s^{3/2}} > 0, \ph t > 0 \ph or  \ph \il{0}{t}d(\sigma_2^+(s) + \sigma_2^-(s)) > 0, \ph t > 0,   \phantom{aaaaaaaaaaaaaaaaaaaaa}
$$
$$
\phantom{aaaaaaaaaaaaaaaaaaaaaaaaaaaaaaaaaaaaaaaaaaa} or \ph \il{0}{t}d(\sigma_3^+(s) + \sigma_3^-(s)) > 0, \ph t > 0,
$$
 then $[Re \hskip 2pt  \widehat{{K}}(\lambda) ]'\ne0, \ph \lambda \ne0
 $.
}

Proof.
It follows from the representation (2.3)   that
$$
[Re \hskip 2pt  \widehat{{K}}_1(\lambda) ]'= -\frac{\sqrt{\pi}\lambda}{2}\int\limits_{0}^\infty\frac{1}{s^{3/2}} e^{-\frac{\lambda^2}{4s}} d(\sigma_1^+(s) + \sigma_1^-(s)), \ph \lambda \in \mathbb{R},
$$
provided the integral $\int\limits_{0}^\infty\frac{|d(\sigma_1^+(s) + \sigma_1^-(s))|}{s^{3/2}} ds$ converges. Since for $\lambda \ne 0$ the function \linebreak $1 - e^{-\frac{\lambda^2}{4s}}, \ph s > 0$ nonincreases on $(0,\infty)$ by the mean value theorem for integrals (see [7, p.  869]) from here we obtain
$$
[Re \hskip 2pt\widehat{{K}}_1(\lambda)]'= \frac{\sqrt{\pi}\lambda}{2}\biggl\{ \int\limits_{0}^\infty[1 - e^{-\frac{\lambda^2}{4s}}] \frac{d(\sigma_1^+(s) + \sigma_1^-(s))}{s^{3/2}} - \int\limits_{0}^\infty\frac{d(\sigma_1^+(s) + \sigma_1^-(s))}{s^{3/2}}\biggr\}=
$$
$$
=\frac{\sqrt{\pi}\lambda}{2}\lim\limits_{N\to \infty}\biggl\{\il{0}{N}[1 - e^{-\frac{\lambda^2}{4s}}]\frac{d(\sigma_1^+(s) + \sigma_1^-(s))}{s^{3/2}}
 - \int\limits_{0}^\infty\frac{d(\sigma_1^+(s) + \sigma_1^-(s))}{s^{3/2}}\bigg\}=
$$
$$
=\frac{\sqrt{\pi}\lambda}{2}\lim\limits_{N\to \infty}\biggl\{\il{0}{\alpha(N,\lambda)}\frac{d(\sigma_1^+(s) + \sigma_1^-(s))}{s^{3/2}}
 - \int\limits_0^\infty\frac{d(\sigma_1^+(s) + \sigma_1^-(s))}{s^{3/2}}\bigg\}=
$$
$$
=-\frac{\sqrt{\pi}\lambda}{2}\lim\limits_{N\to \infty}\il{\alpha(N,\lambda)}{\infty} \frac{d(\sigma_1^+(s) + \sigma_1^-(s))}{s^{3/2}},
$$
where $\alpha(N,\lambda) \in (0,N], \ph N> 0, \ph \lambda \in \mathbb{R}$. By the condition  A) from here it follows that $Re \hskip 2pt \widehat{\mathcal{K}}(\lambda)$ non decreases on $(-\infty,0]$ and non increases on $[0,\infty)$.
Let us show that $Re\hskip 2pt \widehat{\mathcal{K}_2(\lambda)}$ non decreases on $(-\infty,0]$ and non increases  on $[0,\infty)$. We have
$$
Re\hskip 2pt [\widehat{{K}_2(\lambda)} ]'= -\int\limits_{0}^\infty e^{-s|\lambda|}\frac{d(\sigma_2^+(s) + \sigma_2^-(s))}{8\pi}
, \phh \lambda > 0. \eqno (2.6)
$$
Since $e^{-s|\lambda|} > 0$ and non increases  by $s$ on $[0,\infty)$ for every $\lambda > 0$ by mean value theorem for integrals we get
$$
\int\limits_{0}^\infty e^{-s|\lambda|}d[\sigma_2^+(s) + \sigma_2^-(s)] = \lim\limits_{N\to \infty}\il{0}{N}e^{-s|\lambda|}d[\sigma_2^+(s) + \sigma_2^-(s)] = \lim\limits_{N\to \infty}\il{0}{\gamma(N,\lambda)}d[\sigma_2^+(s) + \sigma_2^-(s)],
$$
where $0< \gamma(N,\lambda) \le N$. As for as $\widehat{{K}_2(\lambda)}$ is a even function From here, (2.6)   and the  conditions B) \ph ($j=2$)  it follows that
$Re\hskip 2pt \widehat{{K}_2(\lambda)}$ non decreases  on $(-\infty,0]$ and non increases   on $[0,\infty)$. By similar way using the inequality B) for   $j=3$ one can show that $Re\hskip 2pt \widehat{{K}_3(\lambda)}$ non decreases  on $(-\infty,0]$ and non increases on $[0,\infty)$.
Therefore, it follows from (2.2)  that  $Re \hskip 2pt \widehat{{K}}(\lambda)$ non decreases  on $(-\infty,0]$ and non increases  on $[0,\infty)$. Obviously if
 $\int\limits_{t}^\infty\frac{d(\sigma_1^+(s) + \sigma_1^-(s))}{s^{3/2}} > 0, \ph t \ge 0$
  \Big($\il{0}{t}d(\sigma_2^+(s) + \sigma_2^-(s)) > 0, \ph t \ge 0, \ph  \il{0}{t}d(\sigma_3^+(s) + \sigma_3^-(s)) > 0, \linebreak t \ge~ 0\Big)$, then $[Re \hskip 2pt  \widehat{{K}_1}(\lambda) ]'\ne0, \ph \lambda \ne0, \ph \Big([Re \hskip 2pt  \widehat{{K}_2}(\lambda) ]'\ne0, \ph \lambda \ne0, \ph   [Re \hskip 2pt  \widehat{{K}_3}(\lambda) ]'\ne0, \linebreak \lambda \ne0\Big)$. Therefore, $[Re \hskip 2pt  \widehat{{K}}(\lambda) ]'\ne0, \ph \lambda \ne0$.
The lemma is proved.

\vskip 10 pt
{\bf Remark 2.1.} {\it The next example shows that from the conditions (1.2), $A)$  and $B)$ does not follow the inequality $K(t)+K(-t)\ge 0, \ph t \ge 0$.
}

\vskip 10pt

{\bf Example 2.1.} {\it Let $K(-t) = K(t)=K_1(t)=\frac{1}{3\sqrt{\pi}}(-e^{-t^2}+8e^{-4t^2})$.
 We have
$$
d(\sigma_1^-(s)+\sigma_1^+(s))=\frac{1}{3\sqrt{\pi}}
d[- \hskip 3pt \theta(s-1)+8 \hskip 3pt \theta(s-4)], \ph
 d(\sigma_j^-(s)+\sigma_j^+(s))\equiv 0, \ph j=2,3.
$$
It is easy to check that the conditions $A)$ and $B)$ are satisfied.
Then using the formula
$$
\int\limits_0^\infty e^{-st^2}dt = \frac{1}{2}s^{-1/2}\sqrt{\pi}
$$
one can verify that the condition (1.2) also is satisfied. Moreover, one can show that
$$
\int\limits_{-\infty}^\infty|K(t)| dt \approx 1.1099  \ph (according \ph  to \ph  google  \ph AI).
$$
 Obviously, in this particular case $K(-t)+K(t)=2K(t)<0$ for all enough large $|t|$. Therefore, the condition $K(t)+K(-t)\ge0, \ph t \ge 0$ in this case is not fulfilled.
Moreover, using the formulae
$$
\int\limits_0^\infty t^2 e^{-st^2}dt = \frac{1}{4}s^{-3/2}\sqrt{\pi}, \phh
 \int\limits_0^\infty t^4 e^{-s t^2}dt=\frac{3\sqrt{\pi}}{8 s^{5/2}}, \ph s >0
$$
one can show that
$$
\int\limits_0^\infty t^2K(t)d t = 0,  \phh \int\limits_0^\infty t^4 K(t) d t=-\frac{9}{32}.
$$
}

Using   (2.2)-(2.5)          one can easily prove the following lemmas

\vskip 10pt

{\bf Lemma 2.2.} {\it Let the following conditions be satisfied
$$
A_\pm) \ph \pm d[\sigma_1j^+(s)-\sigma_j^-(s)]\ge0, \ph s\ge 0, \ph j=1,2,3,\ph \pm\sum\limits_{j=1}^3 d[\sigma_j^+(s)-\sigma_j^-(s)]>0.
$$
Then
$$
\pm Im \hskip 2pt \widehat{K}(\lambda) > 0, \ph \lambda > 0.
$$
}

\vskip 10pt

{\bf Lemma 2.3.} {\it If
$$
d[\sigma_j^+(s)-\sigma_j^-(s)]\equiv 0, \ph s\ge 0, \ph j=1,2,3,
$$
then $Im \hskip 2pt \widehat{K}(\lambda)\equiv 0$.
}

Let $\mathcal{A}(\lambda)\equiv\mu-\int\limits_{-\infty}^\infty e^{i\lambda t} K(t) d t, \ph -\infty\le\lambda\le\infty$ and let
 the conditions $A)$ and $B)$ hold.

 Then, obviously, by Lemma 2.1 in the "dissipative" \ph case ($\mu>1$) we have $Re \hskip 2pt \mathcal{A}(\lambda) >~ 0, \linebreak  -\infty \le \lambda \le \infty$. Therefore in this case
$$
\mathcal{A}(\lambda) \ne 0, \ph -\infty \le \lambda \le \infty, \phh ind _{-\infty \le \lambda \le \infty}\mathcal{A}(\lambda) = 0.
$$
In the  "conservative"  \ph  case  ($\mu =1$) the function  $\mathcal{A}(\lambda)$ degenerates at $\lambda = 0$.

Let $\nu_1(K)$ converges absolutely and $\nu_1(K) \ne 0$. Then by (2.1) and  Theorem 2.2  if the conditins  of Lemma 2.1   hold, then
$$
\mathcal{A}(\lambda) = \frac{\lambda}{\lambda+ i} \mathcal{B}(\lambda), \ph \mathcal{B}(\lambda) \in \mathbb{W}, \ph \mathcal{B}(\lambda) \ne 0, \ph  -\infty \le \lambda \le \infty. \eqno   (2.7)
$$
We show that if $\nu_1(K) > 0,$ then
$$
ind_{-\infty \le \lambda \le \infty} \mathcal{B}(\lambda) = 0. \eqno  (2.8)
$$
By     (2.7)  we have $\mathcal{B}(\lambda) = \mathcal{A}(\lambda) + \frac{i}{\lambda} \mathcal{A}(\lambda), \ph \lambda \ne 0$. Then
$$
\mathcal{B}(0) = \mathcal{A}(0)  + \lim\limits_{\lambda\to 0} \frac{i}{\lambda} \mathcal{A}(\lambda) = \nu_1(K) > 0 . \eqno (2.9)
$$
By Lemma 2.1 we have $-\frac{\pi}{2} < arg \mathcal{A}(\lambda) < \frac{\pi}{2}, \ph \lambda \ne 0.$
 This inequality is also true for $arg \frac{i+ \lambda}{\lambda}$. Hence, by (2.7)  $-\pi < arg \mathcal{B}(\lambda) < \pi, \ph \lambda \ne 0$. This together with (2.9) implies (2.8).
Let us show that if $\nu_1(K) < 0,$ then
$$
ind_{-\infty \le \lambda \le \infty} \mathcal{B}(\lambda) = -1. \eqno (2.10)
$$
Consider the function $\widetilde{\mathcal{A}}(\lambda) \equiv 1 - \int\limits_{-\infty}^{\infty}e^{i\lambda t}\widetilde{K}(t)d t, \ph -\infty \le \lambda \le \infty,$ where $\widetilde{K}(t) = K(-t)$. Obviously, $\nu_1(\widetilde{K}) = -\nu_1(k) > 0$. Then by already proven
$$
\widetilde{\mathcal{A}}(\lambda) = \frac{\lambda}{\lambda+ i} \widetilde{\mathcal{B}}(\lambda), \ph \widetilde{\mathcal{B}}(\lambda) \in \mathbb{W}, \ph \widetilde{\mathcal{B}}(\lambda) \ne 0, \ph  -\infty \le \lambda \le \infty \eqno (2.11)
$$
and
$$
ind_{-\infty \le \lambda \le \infty} \widetilde{\mathcal{B}}(\lambda) = 0.  \eqno (2.12)
$$
Obviously, $\widetilde{\mathcal{A}}(\lambda) = \mathcal{A}(-\lambda), \ph -\infty \le \lambda \le \infty.$ Hence, by (2.7)  and  (2.11)  we have
$$
\mathcal{B}(-\lambda) = \Bigl(\frac{\lambda - i}{\lambda + i}\Bigr)\widetilde{\mathcal{B}}(\lambda), \ph  -\infty \le \lambda \le \infty.
$$
This together with (2.12)   implies (2.10). Let us cosider the case $\mu=1, \ph \nu_1(K)=0, \linebreak \nu_2(K)\ne~ 0$. Then by (2.1) and Theorem 2.2
$$
\mathcal{A}(\lambda)=\Big(\frac{\lambda}{\lambda+i}\Big)^2\mathcal{B}
(\lambda), \ph \mathcal{B}(\lambda)\in\mathbb{W}.
$$
From here  we get
$$
\mathcal{A}(\lambda)=\frac{\lambda^2}{\lambda^2+1}\mathcal{B}_1
(\lambda), \ph \mathcal{B}_1(\lambda)\in\mathbb{W}. \eqno (2.13)
$$
where $\mathcal{B}_1(\lambda)=\frac{\lambda-i}{\lambda+i}\mathcal{B}(\lambda), \ph -\infty\le\lambda\le\infty$.

If the conditions of Lemma 2.1 hold, then  $Re \mathcal{B}_1(\lambda) > 0, \ph \lambda \ne 0$. Let us show that $\mathcal{B}_1(0)\ne 0$.
We have $\mathcal{B}_1(0)=\lim\limits_{\lambda\to 0}\frac{\lambda^2+1}{\lambda^2}\mathcal{A}(\lambda)$. Using the L'hospital rule from here we obtain
$
\mathcal{B}_1(0)=\frac{1}{2}\nu_2(K)\ne 0.
$
Therefore,
$$
\mathcal{B}_1(\lambda)\ne0, \ph -\infty\le \lambda\le\infty. \eqno (2.14)
$$
From here and (2.13) we get
$$
\mathcal{A}(\lambda)=\frac{\lambda^2}{\lambda^2+1}\mathcal{B}_1
(\lambda), \ph \mathcal{B}_1(\lambda)\in\mathbb{W}, \ph
 \mathcal{B}_1(\lambda)\ne 0, \ph -\infty\le\lambda\le\infty.
$$
where $\mathcal{B}_1(\lambda)=\frac{\lambda-i}{\lambda+i}\mathcal{B}(\lambda), \ph -\infty\le\lambda\le\infty$. Let us show that
$$
ind_{-\infty \le \lambda \le \infty}\mathcal{B}_1(\lambda) = 0. \eqno(2.15)
$$
 We have
$$
Re \mathcal{B}_1(\lambda) = \frac{\lambda^2+1}{\lambda^2}Re \mathcal{A}(\lambda) \ge 0, \ph -\infty\le\lambda\le\infty.
$$
Therefore,  (2.15)  is valid. From (2.15)   we get
$$
ind_{-\infty \le \lambda \le \infty}\mathcal{B}(\lambda) = -1.
$$
\vskip 10 pt

{\bf Remark 2.2} {\it From the proof of  (2.14)  is seen that $\nu_2(K)>0$. Moreover, from the conditions of Lemma 2.1 is seen that in the ``conservative'' case
$\nu_2(K)\ge 0$ and the equality $\nu_2(K)=0$ may  holds  (see above Example 2.1).
}

Let $\nu_1(K), \nu_2(K)$ and $\nu_3(K)$ converge absolutely. We consider the case $\mu=1, \ph \nu_1(K)=\nu_2(K)=0, \ph \nu_3(K)\ne 0$. Asuume the conditions of Lemma 2.1 are satisfied. Then by  (2.1) and Theorem 2.2
$$
\mathcal{A}(\lambda)=\Big(\frac{\lambda}{\lambda+i}\Big)^3\mathcal{B}(\lambda), \ph \mathcal{B}(\lambda)\in\mathbb{W}, \ph \mathcal{B}(\lambda)\ne 0, \ph \lambda\in[-\infty,0)\cup(0,\infty]. \eqno (2.16)
$$
Let us show that $\mathcal{B}(0)\ne 0$. We have
$$
\mathcal{B}(0)=\lim\limits_{\lambda \to 0}\Big(1+\frac{3i}{\lambda}-\frac{3}{\lambda^2}-\frac{i}{\lambda^3}\Big)\mathcal{A}(\lambda).
$$
Using the L'hospitale's rule from here we get $\mathcal{B}(0)=-\frac{\nu_3(K)}{6}\ne 0$. This together with (2.16)     imples that
$$
\mathcal{A}(\lambda)=\Big(\frac{\lambda}{\lambda+i}\Big)^3\mathcal{B}(\lambda), \ph \mathcal{B}(\lambda)\in\mathbb{W}, \ph \mathcal{B}(\lambda)\ne 0, \ph -\infty\le\lambda\le\infty. \eqno (2.17)
$$
From here we obtain
$$
\mathcal{A}(\lambda)=\frac{\lambda}{\lambda+i}\frac{\lambda^2}{\lambda^2+1}\mathcal{B}_1(\lambda), \ph -\infty\le \lambda\le \infty, \eqno (2.18)
$$
where $\mathcal{B}_1(\lambda)=\frac{\lambda-i}{\lambda+i}\mathcal{B}(\lambda), \ph -\infty\le\lambda\le\infty$. Let us show that if $\nu_3(K)<0$, then
$$
ind_{-\infty\le\lambda\le\infty}\mathcal{B}_1(\lambda)=0. \eqno (2.19)
$$
From (2.18)   we have
$$
\mathcal{B}_1(\lambda)=\Big(1+\frac{i}{\lambda}\Big)\Big(1+\frac{1}{\lambda^2}\Big)\mathcal{A}(\lambda), \ph\lambda\in[-\infty,0)\cup(0,\infty]. \eqno (2.20)
$$
Obviously, $-\frac{\pi}{2}< arg (1+\frac{i}{\lambda})<\frac{\pi}{2}, \ph \lambda\in[-\infty,0)\cup(0,\infty]$ and (due to the conditions of Lemma 2.1  ) $-\frac{\pi}{2}< arg\Big(1+\frac{1}{\lambda^2}\Big)\mathcal{A}(\lambda)<\frac{\pi}{2}, \ph \lambda\in[-\infty,0)\cup(0,\infty]$ . This together with (2.20)  implies that $-\pi< arg\mathcal{B}_1(\lambda)<\pi, \ph \lambda\in[-\infty,0)\cup(0,\infty]$.
Then sinse $\mathcal{B}_1(0)=-\mathcal{B}(0)=-\frac{\nu_3(K)}{6}>0$ we
conclude that  $-\pi< arg\mathcal{B}_1(\lambda)<\pi, \ph -\infty\le\lambda\le\infty$. Therefore, (2.19)   is valid. Let us show that if $\nu_3(K)>0$, them
$$
ind_{-\infty\le\lambda\le\infty}\mathcal{B}_1(\lambda)=-1. \eqno (2.21)
$$

Consider the function $\widetilde{\mathcal{A}}(\lambda) \equiv 1 - \int\limits_{-\infty}^{\infty}e^{i\lambda t}\widetilde{K}(t)d t, \ph -\infty \le \lambda \le \infty,$ where $\widetilde{K}(t) = K(-t)$. Obviously, $\nu_3(\widetilde{K}) = -\nu_3(k) <0$. Then by already proven
$$
\widetilde{\mathcal{A}}(\lambda) = \frac{\lambda}{\lambda+ i}\frac{\lambda^2}{1+\lambda^2} \widetilde{\mathcal{B}_1}(\lambda), \ph \widetilde{\mathcal{B}_1}(\lambda) \in \mathbb{W}, \ph \widetilde{\mathcal{B}_1}(\lambda) \ne 0, \ph  -\infty \le \lambda \le \infty \eqno (2.22)
$$
and
$$
ind_{-\infty \le \lambda \le \infty} \widetilde{\mathcal{B}_1}(\lambda) = 0.  \eqno (2.23)
$$
Obviously, $\widetilde{\mathcal{A}}(\lambda) = \mathcal{A}(-\lambda), \ph -\infty \le \lambda \le \infty.$ Hence,  by (2.22)   we have
$$
\mathcal{B}_1(-\lambda) = \Bigl(\frac{\lambda - i}{\lambda + i}\Bigr)\widetilde{\mathcal{B}_1}(\lambda), \ph  -\infty \le \lambda \le \infty.
$$
This together with (2.23)  implies  (2.21). From (2.17),        (2.19) and (2.21)  we derive
$$
ind_{-\infty \le \lambda \le \infty} \mathcal{B}(\lambda) = \sist{-1, \ph if \ph \nu_3(K)<0,}{-2, \ph if \ph \nu_3(K)>0.}
$$
Assume $\nu_m(K),  \ph m=1,2,3,4$ converge absolutely. We consider the case $\mu=1, \ph \nu_j(K)=0, \ph j=1,2,3, \ph \nu_4(K)\ne0$. Assume the conditions of Lemma 2.1 are satisfied. Then by  (2.1) and Theorem 2.2
$$
\mathcal{A}(\lambda)=\Big(\frac{\lambda}{\lambda+i}\Big)^4\mathcal{B}(\lambda), \ph \mathcal{B}(\lambda)\in\mathbb{W}, \ph \mathcal{B}(\lambda)\ne 0, \ph \lambda\in[-\infty,0)\cup(0,\infty]. \eqno (2.24)
$$
Let us show that $\mathcal{B}(0)\ne 0$. We have
$$
\mathcal{B}(0)=\lim\limits_{\lambda \to 0}\Big(1+\frac{4i}{\lambda}-\frac{6}{\lambda^2}-\frac{3i}{\lambda^3}+\frac{1}{\lambda^4}\Big)\mathcal{A}(\lambda).
$$
Using the L'hospitale's rule from here we get $\mathcal{B}(0)=-\frac{\nu_4(K)}{24}\ne 0$. This together with (2.24)   imples that
$$
\mathcal{A}(\lambda)=\Big(\frac{\lambda}{\lambda+i}\Big)^4\mathcal{B}(\lambda), \ph \mathcal{B}(\lambda)\in\mathbb{W}, \ph \mathcal{B}(\lambda)\ne 0, \ph -\infty\le\lambda\le\infty.
$$
From here we obtain
$$
\mathcal{A}(\lambda)=\frac{\lambda^4}{(\lambda^2+1)^2}\mathcal{B}_2(\lambda), \ph -\infty\le \lambda\le \infty, \eqno (2.25)
$$
where $\mathcal{B}_2(\lambda)=\Big(\frac{\lambda-i}{\lambda+i}\Big)^2\mathcal{B}(\lambda), \ph -\infty\le\lambda\le\infty$. Let us show that
$$
ind_{-\infty\le\lambda\le\infty}\mathcal{B}_2(\lambda)=0. \eqno (2.26)
$$
Since under the conditions of Lemma 2.1 $-\frac{\pi}{2}\le arg \frac{(\lambda^2+1)^2}{\lambda^4}\mathcal{A}(\lambda)\le \frac{\pi}{2}, \ph -\infty\le\lambda\le\infty$ from (2.25) is seen that $-\frac{\pi}{2}\le arg\mathcal{B}_2(\lambda)\le \frac{\pi}{2}, \ph -\infty\le\lambda\le\infty$.  This together with (2.25)  implies (2.26). From (2.26) we immediately get
$$
ind_{-\infty\le\lambda\le\infty}\mathcal{B}(\lambda)=-2.
$$
\vskip 10 pt

{\bf Remark 2.3} {\it From the proof of  (2.34) is seen that $\nu_4(K)<0$.
}



Let     $0<\mu<1$ (the ``supercritical''  case).
Then if the conditins  $A_+)$
of Lemma 2.2 hold,  then the courve $\mathcal{A}(\lambda)$ connects the point $\mu+i0$  with $\mu-1+i 0$ when $\lambda$ varies from $-\infty$ to $0$, remaining in the uper half plane of $\mathbb{C}$. Therefore,  when $\lambda$ varies from $-\infty$ to $\infty$ the courve $\mathcal{A}(\lambda)$ movies around the origin of the complex plane in the contrclockvise direction and returms to the starting point $\mu+i0$.
  Therefore, in this case
$$
\mathcal{A}(\lambda) \ne 0, \ph -\infty \le \lambda \le \infty, \phh ind _{-\infty \le \lambda \le \infty}\mathcal{A}(\lambda) =1.
$$
If the conditions  $A_-)$
of Lemma 2.2 hold,  then the courve $\mathcal{A}(\lambda)$ starting from the point $\mu+i0$ movies around the origin of the complex plane in the clockvise direction and returms to the starting point $\mu+i0$ when $\lambda$ movies from $-\infty$ to $\infty$.  Therefore, in this case
$$
\mathcal{A}(\lambda) \ne 0, \ph -\infty \le \lambda \le \infty, \phh ind _{-\infty \le \lambda \le \infty}\mathcal{A}(\lambda) =-1 .
$$
Obviously, if the conditions of Lemmas 2.1 and 2.3 hold, then  $\mathcal{A}(\pm\alpha)=0$ for some $\alpha > 0$ and other zeroes it has not.  Let us show that
$$
(R_\alpha K)(t)\in L^1
$$
Let $t>0$ then $(R_\alpha K_1)(t)=ie^{-i\alpha t}\int\limits_t^\infty e^{i\alpha s} K_1(s) ds=ie^{-i\alpha t}\int\limits_t^\infty e^{i\alpha s} d s \int\limits_0^\infty e^{-\xi s^2}d\sigma_1^+(\xi)= \linebreak
=\int\limits_0^\infty e^{-\xi s^2}d\sigma_1^+(\xi)\int\limits_0^\xi e^{i\alpha s} d s =i /\alpha \int\limits_0^\infty  (1- e^{\alpha\xi})  e^{-\xi s^2}d\sigma_1^+(\xi)$. Since $  (1- e^{\alpha\xi})$ is bounded it follows from here that $(R_\alpha K_1)(t)$ is absolutely integrable on $[0,\infty)$. By amalogy can be proved that  $(R_\alpha K_1)(t)$ is absolutely integrable on $(-\infty,0]$. Threefore,  $(R_\alpha K_1)(t) \in L^1$. By analogy can be proved that  $(R_\alpha K_j)(t) \in L^1, \ph j=1,2$. Therefore,  $(R_\alpha K)(t) \in L^1$.
 By  Theorem 2.2 it follows  from here that
$$
\mathcal{A}(\lambda)=\frac{\lambda^2-\alpha^2}{[\lambda+i]^2}
\mathcal{B}(\lambda), \ph \mathcal{B}(\lambda)\in\mathbb{W}. \eqno  (2.27)
$$
If the conditions of Lemma 2.1 hold, then
$$
\mathcal{B}(\lambda)\ne 0, \ph -\infty\le\lambda\le\infty.
$$
Let us show that
$$
ind _{-\infty \le \lambda \le \infty}\mathcal{B}(\lambda) =-1. \eqno (2.28)
$$
We set $\mathcal{B}_1(\lambda)\equiv \bigl(\frac{\lambda -i}{\lambda+i}\bigr) \mathcal{B}(\lambda), \ph -\infty\le \lambda\le\infty$.  Then it follows from (2.27)   that
$$
\mathcal{A}(\lambda)=\frac{\lambda^2-\alpha^2}{\lambda^2+1}
\mathcal{B}_1(\lambda), \ph \mathcal{B}_1(\lambda)\in\mathbb{W}, \ph
 \mathcal{B}_1(\lambda)\ne 0, \ph -\infty
 \le \lambda \le \infty. \eqno (2.29)
$$
Under the comditions of Lemma 2.3  we have $Im \hskip 2pt \mathcal{A}(\lambda)=0, \ph -\infty\le\lambda\le\infty$. This together with (2.29)  implies that $Im \hskip 2pt \mathcal{B}_1(\lambda)=0, \ph -\infty\le\lambda\le\infty$. Therefore, \linebreak $ind _{\infty \le \lambda \le \infty}\mathcal{B}_1(\lambda) =0$. From here, (2.27) and (2.29)   follows  (2.28).

 Let $\mu<0$ (the ``antidisipative'' case). Then then under the restrictions $A)$ and $B)$ $Re \mathcal{A}(\lambda) \le\mu$. Therefore, in this case
$$
\mathcal{A}(\lambda) \ne 0, \ph -\infty \le \lambda \le \infty, \phh ind _{-\infty \le \lambda \le \infty}\mathcal{A}(\lambda) = 0.
$$
 Summarizing the obtained we come to the following theorem.
\vskip 10pt

{\bf Theorem 2.3.} {\it  The following statements are valid.

\noindent
I) If $\mu\in(-\infty,0)\cup(1,\infty)$ and the conditions    $A)$   and $B)$ hold, then
$$
\mathcal{A}(\lambda) \ne 0, \ph -\infty \le \lambda \le \infty, \phh ind _{-\infty \le \lambda \le \infty}\mathcal{A}(\lambda) = 0.
$$

\noindent
II) If $\mu=1$, the conditions    $A), \ph B)$ and $C)$   hold, $\nu_1(K)$ converges absolutely and $\nu_1(K)\ne0$, then
$$
\mathcal{A}(\lambda) = \frac{\lambda}{\lambda+ i} \mathcal{B}(\lambda), \ph \mathcal{B}(\lambda) \in \mathbb{W}, \ph \mathcal{B}(\lambda) \ne 0, \ph  -\infty \le \lambda \le \infty, \eqno (2.30)
$$
$$
ind _{-\infty \le \lambda \le \infty}\mathcal{B}(\lambda) =\sist{-1, \ph if \ph \nu_1(K)<0,}{0, \ph if \ph  \nu_1(K)>0.} \eqno (2.31)
$$
\noindent
III) If $\mu=1$,  the conditions    $A), \ph B)$ and $C)$   hold, $\nu_1(K)$ and $\nu_2(K)$ converge\ absolutely, \linebreak $\nu_1(K)=0$ and $\nu_2(K)\ne0$, then
$$
\mathcal{A}(\lambda) = \frac{\lambda^2}{(\lambda+ i)^2} \mathcal{B}(\lambda), \ph \mathcal{B}(\lambda) \in \mathbb{W}, \ph \mathcal{B}(\lambda) \ne 0, \ph  -\infty \le \lambda \le \infty, \eqno (2.32)
$$
$$
ind_{-\infty \le \lambda \le \infty}\mathcal{B}(\lambda) = -1. \eqno (2.33)
$$
\noindent
IV) If  $\mu=1$,  the conditions    $A), \ph B)$ and $C)$   hold,  $\nu_j(K), \ph j=1,2,3$ converge  absolutely, \linebreak $\nu_1(K)=\nu_2(K)=0$ and $\nu_3(K)\ne0$, then
$$
\mathcal{A}(\lambda) = \frac{\lambda^3}{(\lambda+ i)^3} \mathcal{B}(\lambda), \ph \mathcal{B}(\lambda) \in \mathbb{W}, \ph \mathcal{B}(\lambda) \ne 0, \ph  -\infty \le \lambda \le \infty,
$$
$$
ind _{-\infty \le \lambda \le \infty}\mathcal{B}(\lambda) =\sist{-1, \ph if \ph \nu_3(K)<0,}{-2, \ph if \ph  \nu_3(K)>0.}
$$
\noindent
V) If  $\mu=1$,  the conditions    $A), \ph B)$ and $C)$   hold,  $\nu_j(K), \ph j=1,2,3,4$ converge absolutely, \linebreak $\nu_1(K)=\nu_2(K)=\nu_3(K)=0$ and $\nu_4(K)\ne0$, then
$$
\mathcal{A}(\lambda) = \frac{\lambda^4}{(\lambda+ i)^4} \mathcal{B}(\lambda),
\ph \mathcal{B}(\lambda) \in \mathbb{W},
\ph \mathcal{B}(\lambda) \ne 0, \ph  -\infty \le \lambda \le \infty,
$$
$$
ind _{-\infty \le \lambda \le \infty}\mathcal{B}(\lambda) =-2.
$$
\noindent
VI) If $0<\mu<1$ amd the condition $A_{\pm})$
holds, then
$$
\mathcal{A}(\lambda) \ne 0, \ph -\infty \le \lambda \le \infty, \phh ind _{-\infty \le \lambda \le \infty}\mathcal{A}(\lambda) =\pm 1.
$$
\noindent
VII) If $d[\sigma_j^+(s)-\sigma_j^-(s)]\equiv 0, \ph s\ge 0, \ph j=1,2,3$, then
$$
\mathcal{A}(\lambda)=\frac{\lambda^2-\alpha^2}{[\lambda+i]^2}
\mathcal{B}(\lambda), \ph \alpha >0, \ph  \mathcal{B}(\lambda)\in\mathbb{W}, \ph \mathcal{B}(\lambda)\ne 0, \ph -\infty\le\lambda\le\infty,
$$
$$
ind _{-\infty \le \lambda \le \infty}\mathcal{B}(\lambda) =-1.
$$
}

\vskip 20pt

{\bf 3. Main reaults.}  Denote by $\mathbb{A}_\mu$ the operator, corresponding to Eq (1.1).  First we consider the  ``disipative''  ($\mu>1$) and the ``antidisipative'' ($\mu<0$) cases. If the comditions of statement I) of Theorem 2.3  are satisfied, then operator $\mathbb{A}_\mu: E_+\to E_+$ is invertible. From here we immediately derive

\vskip 10pt

{\bf Theorem 3.1.} {\it Let $\mu\in(-\infty,0)\cup(1,\infty)$ and the conditions $A), \ph B)$ hold. Then for every $f(t)\in E_+$ Eq. (1.1) has a solution $\phi(t)\in E_+$. The corresponding homogeneous equation has only a trivial solution in $E_+$.
}

\vskip 10pt

{\bf Remark 3.1.} {\it The next example shows, that in general from the conditions $|\mu|>~1, \linebreak  A), \ph B)$ and (1.2) does not follow that  the operator $\frac{1}{\mu}\mathbb{A}_\mu - I$ is a contracting mapping in $L^p_+, \ph 1\le p\le\infty$ (see also above, Example 2.1).
}

\vskip 10pt

{\bf Example 3.1.} {\it Let $K(t)=\frac{2000}{\sqrt{\pi}}(-e^{-t^2}+2.001 e^{-4 t^2})$. Then it is not difficult to verify that the conditions (1.2) $A)$ and $B)$ are satisfied, and
$$
\int\limits_{-\infty}^\infty|K(t)|d t \approx 1291.35 \ph (according \ph to \ph  google \ph AI).
$$
}

Now let us explore the ``conservative'' ($\mu=1$) case.
Assume the conditions of  statement II)  of Theorem 2.3 hold.
 Then by Theorems 2.1 and 2.3  for every $f(t)\in E_+$ Eq. (1.1) has a solution  $\phi(t) \in \widetilde{E}_+(\frac{\lambda}{\lambda+i})$. The corresponding homogeneous equation has: the unique trivial solution in   $\widetilde{E}_+(\frac{\lambda}{\lambda+i})$  if $\nu_1(K)>0$, and  up to an arbitrary multiplier  the unique nontrivial solution  \linebreak  $\phi_0(t) \in~ \widetilde{E}_+(\frac{\lambda}{\lambda+i})$  , if $\nu_1(K)<0$. It follows from (2.30) and (2.31) that
$$
\mathcal{A}(\lambda) = \frac{\lambda}{\lambda- i} \widetilde{\mathcal{B}}(\lambda), \ph \widetilde{\mathcal{B}}(\lambda) \in \mathbb{W}, \ph \widetilde{\mathcal{B}}(\lambda) \ne 0, \ph  -\infty \le \lambda \le \infty,
$$
$$
ind _{-\infty \le \lambda \le \infty}\widetilde{\mathcal{B}}(\lambda) =\sist{0, \ph if \ph \nu_1(K)<0,}{1, \ph if \ph  \nu_1(K)>0.}
$$
By Theorem 2.1 it follows from here that, if $\nu_1(K)<0$, then $\phi_0(t)\in C_+$ and for every $f(t)\in \overline{E}_+(\frac{\lambda}{\lambda-i})$ Eq. (1.1) has a solution in $E_+$. Moreover if $\nu_1(K)>0$, then for every $f(t)\in \overline{M}_+(\frac{\lambda}{\lambda-i}) \ph( f(t)\in \overline{M}_+(\frac{\lambda}{\lambda-i})$) Eq. (1.1) has a solution in $M_+ \ph (C_+)$.
Thus we have proved the following theorem.

\vskip 10pt

{\bf Theorem 3.2.} {\it Let $\mu=1$, the conditions    $A), \ph B)$ and $C)$   hold, $\nu_1(K)$ converges absolutely and $\nu_1(K)\ne0$. Then for every $f(t)\in E_+$
 Eq. (1.1) has a solution $\phi(t)\in\widetilde{E}_+(\frac{\lambda}{\lambda+i})$.  If $\nu_1(K)>0$, then the corresponding homogeneous equation has: the unique trivial solution in   $\widetilde{E}_+(\frac{\lambda}{\lambda+i})$. If $\nu_1(K)<0$, then Eq. (1.1) has up to an arbitrary multiplier  the unique nontrivial solution    $\phi_0(t) \in~C_+$. Moreover, if $\nu_1(K)>0$, then for every $f(t)\in \overline{M}_+(\frac{\lambda}{\lambda-i}) \ph( f(t)\in \overline{M}_+(\frac{\lambda}{\lambda-i})$) Eq. (1.1) has a solution in $M_+ \ph (C_+)$.
}

Assume the conditions of  statement III)  of Theorem 2.3 hold.
Then by Theorems 2.1 amd 2.3  for every $f(t)\in E_+$ Eq. (1.1) has a solution  $\phi(t) \in \widetilde{E}_+(\frac{\lambda^2}{(\lambda+i)^2})$. The corresponding homogeneous equation has up to an arbitrary multiplier the unique solution $\phi_0(t)\in\widetilde{E}_+(\frac{\lambda^2}{(\lambda+i)^2})$.  It follows from (2.32) and (2.33) that
$$
\mathcal{A}(\lambda) = \frac{\lambda}{\lambda+ i} \widetilde{\mathcal{B}}(\lambda)\frac{\lambda}{\lambda-i}, \ph \widetilde{\mathcal{B}}(\lambda) \in \mathbb{W}, \ph \widetilde{\mathcal{B}}(\lambda) \ne 0, \ph  -\infty \le \lambda \le \infty,
$$
$$
ind_{-\infty \le \lambda \le \infty}\widetilde{\mathcal{B}}(\lambda) =0.
$$
By Theorem 2.1 it follows from here that, if $f(t)\in \overline{E_+}(\frac{\lambda}{\lambda-i})$,  then $\phi(t)\in \widetilde{E}_+(\frac{\lambda}{\lambda+ i})$ and $\phi_0(t)\in\widetilde{C}_+(\frac{\lambda}{\lambda+i})$. The obtained result we formulate in the following form.

\vskip 10 pt

{\bf Theorem 3.3.} {\it Assume   $\mu=1$,  the conditions    $A), \ph B)$ and $C)$   hold, $\nu_1(K)$ and $\nu_2(K)$ converge absolutely,  $\nu_1(K)=0$ and $\nu_2(K)\ne0$. Then  for every $f(t)\in E_+$ Eq. (1.1) has a solution  $\phi(t) \in \widetilde{E}_+(\frac{\lambda^2}{(\lambda+i)^2})$. The corresponding homogeneous equation has up to an arbitrary multiplier the unique solution  $\phi_0(t)\in\widetilde{C}_+(\frac{\lambda}{\lambda+i})$. Moreover, if $f(t)\in \overline{E}_+(\frac{\lambda}{\lambda-i})$,  then $\phi(t)\in \widetilde{E}_+(\frac{\lambda}{\lambda+ i})$.
}

Using statements IV) and V) of Theorem 2.3 on the basis of Theorem 2.1 by analogy with the proofs of Theorems 3.2 and 3.3 one can prove the following theorems.

\vskip 10 pt

{\bf Theorem 3.4.} {\it Assume  $\mu=1$,  the conditions    $A), \ph B)$ and $C)$   hold,  $\nu_j(K), \ph j=1,2,3$ converge  absolutely,  $\nu_1(K)=\nu_2(K)=0$ and $\nu_3(K)\ne0$. Then for every  $f(t)\in E_+$ Eq. (1.1) has a solution  $\phi(t) \in \widetilde{E}_+(\frac{\lambda^3}{(\lambda+i)^3})$.
If $\nu_3(K)<0$, then the corresponding homogeneous equation has  up to an arbitrary multiplier  the unique nontrivial solution    $\phi_0(t) \in~C_+$. If  $\nu_3(K)>0$, then the corresponding homogeneous equation has two linearly independent solutions $\phi_0(t), \phi_1(t)\in C_+$.
 Moreover, if $\nu_1(K)>0$, then for every $f(t)\in \overline{M}_+(\frac{\lambda}{\lambda-i}) \ph( f(t)\in \overline{C}_+(\frac{\lambda}{\lambda-i})$) Eq. (1.1) has a solution in $M_+ \ph (C_+)$.
}

\phantom{aaaaaaaaaaaaaaaaaaaaaaaaaaaaaaaaaaaaaaaaaaaaaaaaaaaaaaaaaaaaaaaaaaaaaaaa}$\blacksquare$

\vskip 10  pt

{\bf Theorem 3.5.} {\it Assume    $\mu=1$,  the conditions    $A), \ph B)$ and $C)$   hold,  $\nu_j(K), \ph j=1,2,3,4$ converge absolutely,  $\nu_1(K)=\nu_2(K)=\nu_3(K)=0$ and $\nu_4(K)\ne0$.
 Then   for every $f(t)\in E_+$ Eq. (1.1) has a solution  $\phi(t) \in \widetilde{E}_+(\frac{\lambda^4}{(\lambda+i)^4})$. The corresponding homogeneous equation has two linearly independent solutions  $\phi_0(t), \phi_1(t)\in  \widetilde{C}_+(\frac{\lambda^2}{(\lambda+i)^2})$.Moreover, if $f(t)\in \overline{E}_+(\frac{\lambda}{\lambda-i})$,  then $\phi(t)\in \widetilde{E}_+(\frac{\lambda^3}{(\lambda- i)^3})$, and if  $f(t)\in \overline{E}_+(\frac{\lambda^2}{(\lambda-i)^2})$,  then $\phi(t)\in \widetilde{E}_+(\frac{\lambda^2}{(\lambda+ i)^2})$.
}

\phantom{aaaaaaaaaaaaaaaaaaaaaaaaaaaaaaaaaaaaaaaaaaaaaaaaaaaaaaaaaaaaaaaaaaaaaaa}$\blacksquare$

Now let us explore the ``supercritical''  ($0<\mu<1$) case. Assume $0<\mu<1$ and the condition $A_+)$ holds. Then by Theorem 2.3,VI) The homogeneous equation, corresponding to Eq. (1.1) has only the trivial solution.  In order that Eq. (1.1) has a solution in $E_+$ it is  necessary and sufficient that $\int\limits_0^\infty f(t) \psi_0(t) d t =0$, where $\psi_0(t)$ is the unique (up to an arbitrary multiplier) solution of the equation
$$
\mu\psi(t)=\int\limits_0^\infty K(\tau -t)\psi(\tau)d\tau, \ph t \ge 0. \eqno (3.1)
$$
If the condition $A_-)$ holds, then according to Theorem 2.3.VI) for every $f(t)\in E_+$ Eq. (1.1) has a solution $\phi(t)\in E_+$, the corresponding  homogeneous  equation has up to an arbitrary multiplier the unique nontrivial solution $\phi_0(t)\in E_+$.
Thas we have proved the following theorem

\vskip 10pt

{\bf Theorem 3.6} {\it Let $0<\mu<1$. If the condition $A_+)$ holds, then the homogeneous equation, corresponding to Eq. (1.1) has only the trivial solution.  In order that Eq. (1.1) has a solution in $E_+$ it is necessary and sufficient that $\int\limits_0^\infty f(t) \psi_0(t) d t =0$, where $\psi_0(t)$ is the unique (up to an arbitrary multiplier) solution of Eq. (3.1).  If the condition $A_-)$ holds, then
for every $f(t)\in E_+$ Eq. (1.1) has a solution $\phi(t)\in E_+$, the corresponding homogeneous   equation has up to an arbitrary multiplier the unique nontrivial solution $\phi_0(t)\in E_+$.
}

Using Theorems 2.1 and 2.3.VII) it easy to prove the following theorem.
\vskip 10pt

{\bf Theorem 3.7.} {\it  Assume $0<\mu<1, \ph d[\sigma_j^+(s)-\sigma_j^-(s)]\equiv 0, \ph s\ge 0, \ph j=1,2,3$. Then there exists $\alpha > 0$  such that for every $f(t)\in E_+$ Eq. (1.1) has a solution
$\phi(t)\in\widetilde{E}_+(\frac{\lambda^2-\alpha^2}{(i+\lambda)^2})$. The corresponding homogeneous equation has up to an arbitrary multiplier the unique solution
$\phi_0(t)\in\widetilde{M}_+(\frac{\lambda\pm\alpha}{\lambda+i})$.
Moreover, if $f(t)\in\overline{E}_+(\frac{\lambda\pm\alpha}{\lambda-i})$, then
$\phi(t)\in\widetilde{E}_+(\frac{\lambda\mp\alpha}{\lambda+i})$.
}

\phantom{aaaaaaaaaaaaaaaaaaaaaaaaaaaaaaaaaaaaaaaaaaaaaaaaaaaaaaaaaaaaaaaaaaaaaaa}$\blacksquare$

 \vskip 10pt

{\bf Remark 3.1.} {\it
The case $\mu=0$ is a special one, reducing Eq. (1.1) to an equation of the first kind, which is a subject for separate study.
}




\vskip 20pt

\centerline{\bf References}

\vskip 20pt

\noindent
1. L. G. Arabadzhyan, The Wiener-Hopf integral equation in the supercritical case.  \linebreak \phantom{a}Vol. 76, 2004, pp.10-17.

\noindent
2. G. A. Grigorian The Wiener-Hopf equation in the supercritical case.
Izv. Nats.  Akad.  \linebreak \phantom{a} Nauk Armenii Mat. 32 (1997), no. 1, 60-74.

\noindent
3. L. G.   Arabadzhyan, N. B. Engibarian, Equations in convolution and nonlinear functional \linebreak \phantom{a} equations (in Russian). Math. Analiz, Itogi Nauki i Tekhnjki, VINITI, Moskow,\linebreak \phantom{a} vol. 22, 1984, pp. 175-244.

\noindent
4. G. A. Grigorian, Solvability of a class of Wiener-Hopf integral equations. Izv. Nats. \linebreak \phantom{a} Akad. Nauk. Armenii. Matematika, vol 21, No. 2, 1996, pp. 21-32.

\noindent
5. S. Presdorf, Some Classes of Singular Equations (North - Holand mathematical library;  \linebreak \phantom{a} v. 17), Amsterdam--New York. Oxford, 1978.

\noindent
6. M. A.  Lavrentev,  B. V. Shabat, The methods of the theory functions of complex  \linebreak \phantom{a}  variables. Moskow, ``Naula'', 1987.

\noindent
7. R. E. Edvards, A Formal Background to Mathematics (Springer-Verlag, New York,  \linebreak \phantom{a}  Heidelberg, Berlin, 1980).  MR0636095

\end{document}